\documentclass[12pt]{amsart}

\newcommand{\Z}{\mathbf Z}
\newcommand{\C}{\mathbf C}

\newcommand{\GG}{GG}

\begin{document}
\title[Non-Gaussian integrals and hypergeometric functions]{Non-Gaussian integrals and general hypergeometric functions}
\author{Alexander Roi Stoyanovsky}
\begin{abstract} 
By a non-Gaussian integral we mean integral of the product of an arbitrary function and exponent of a polynomial. 
We develop a theory of such integrals, which generalizes and simplifies the theory of general hypergeometric 
functions in the sense of I.~M.~Gelfand et~al. 
\end{abstract}
\maketitle

\section*{Introduction}

Let 
\begin{equation}
P(t_1,\ldots,t_n)=\sum\limits_{\omega=(\omega^1,\ldots,\omega^n)}c_\omega t^\omega,\ \ \ t^\omega=t_1^{\omega^1}\ldots t_n^{\omega^n},
\end{equation}
be a polynomial in $n$ variables with complex coefficients, and let $\alpha(t_1$, $\ldots$, $t_n)$ be an arbitrary (possibly multi-valued) function. By the {\it non-Gaussian integral transform of} $\alpha$ or simply by
the {\it non-Gaussian integral} we mean the integral
\begin{equation}
I_\alpha(P)=I_\alpha(c_\omega)=\oint e^{P(t_1,\ldots,t_n)}\alpha(t_1,\ldots,t_n)dt_1\ldots dt_n,
\end{equation}
considered as a function of the coefficients $c_\omega$ of the polynomial $P$. Here the symbol $\oint$ means that integration goes over a real $n$-dimensional oriented contour without boundary
in the $n$-dimensional complex space, equipped with a choice of a continuous branch of the integrated function. 
Since the integration contour can be non-unique, integral (2) is, in general, a multi-valued function. 

In the last decades it became increasingly clear that theory of non-Gaussian integrals plays a fundamental role in mathematics and phys\-ics. 
Some insights to this theory were given in 
[1,  2, 7--12], but the main contribution is due to I.~M.~Gelfand et~al. [3--5]. They have shown that theory of non-Gaussian integrals is closely related with the theory of 
general hypergeometric functions invented by I.~M.~Gelfand et~al. 

In the present paper we generalize and simplify the theory of general hypergeometric functions.
Namely, we study general integrals of the form
\begin{equation}
I_\beta(P_1,\ldots,P_k)=\oint\beta(P_1,\ldots,P_k,t_1,\ldots,t_n)dt_1\ldots dt_n,
\end{equation}
where $P_1,\ldots,P_k$ are polynomials in $t_1,\ldots,t_n$, and $\beta(y_1$, $\ldots$, $y_k$, $t_1$, $\ldots$, $t_n)$ is an arbitrary (possibly multi-valued) function. We call 
these integrals, considered as functions of coefficients of the polynomials 
$P_1,\ldots,P_k$, by {\it \emph(general\emph) hypergeometric functions}. They 
are expressed through non-Gaussian integral (2) as 
\begin{equation}
I_\beta(P_1,\ldots,P_k)=I_\alpha(P),
\end{equation}
where 
\begin{equation}
P(\lambda_1,\ldots,\lambda_k,t_1,\ldots,t_n)=\lambda_1 P_1+\ldots+\lambda_k P_k,
\end{equation} 
$\lambda_1,\ldots,\lambda_k$ are additional variables (the {\it Cayley trick}), and $\alpha(\lambda_1$, $\ldots$, $\lambda_k$, $t_1$, $\ldots$, $t_n)$ is the inverse
Fourier--Laplace transform of $\beta(y_1$, $\ldots$, $y_k$, $t_1$, $\ldots$, $t_n)$ with respect to $y_1,\ldots,y_k$.
The theory of these integrals 
is a generalization and simplification of the theory of general hypergeometric functions due to Gelfand et~al. 
This generalization is similar to passing from toric algebraic geometry to general algebraic geometry. 
The only more general point in the theory of Gelfand et~al. is that in their theory the polynomials $P$ and $P_1,\ldots,P_k$ can be Laurent polynomials, i.~e. polynomials in 
$t_j$ and $t_j^{-1}$.

Let us say a few words about foundations of our theory. In the paper we use, without formalization, the notion ``arbitrary function'' (of several variables). 
The general intuition of function goes back to Euler and, for more modern times, to I.~M. Gelfand. Informally, a function is a dependence of a variable on several variables,
in particular, given by compositions of operations $+,-,\times,:$, derivative, integral, and solving equations.
The most close set theory formalizations of this notion are the notions of  
distribution [6], multi-valued analytical function, $D$-module, and sheaf.

The paper is organized as follows.
In \S1 we give examples of non-Gauss\-ian integrals (2).

In \S2 we give examples of general hypergeometric functions (3). 

In \S3 we derive systems of equations satisfied by the non-Gaussian integrals and by general hypergeometric functions. 
Corollaries of these systems are $A$-hyper\-geo\-metr\-ic systems and $GG$-systems defined and studied in [3--5].

Finally, in \S4 we give power series expansions of non-Gaussian integrals and general hypergeometric functions.

{\bf Acknowledgment.} I would like to thank A. Gemintern for interest to this work.

\section{Examples of non-Gaussian integrals}

Recall that by the non-Gaussian integral we mean integral (2).

{\it Example} 1.1. Let $P=c_1t_1+\ldots+c_nt_n$
 be a linear form in $t_1,\ldots,t_n$, then 
\begin{equation}
I_\alpha(c_1t_1+\ldots+c_nt_n)=\widehat\alpha(c_1,\ldots,c_n)
\end{equation} 
is the Fourier--Laplace transform of $\alpha(t_1,\ldots,t_n)$. 

{\it Example} 1.2. Let $P$ be a quadratic expression in $t_1,\ldots,t_n$, then 
integral (2) is well studied and called the {\it Gaussian integral transform of} $\alpha$ or simply by the {\it Gaussian integral}. For description of the image of the Gaussian integral transform, see [9, 11]
(the results of [9] are announced in [10]). 

{\it Example} 1.3. Let $\alpha(t_1,\ldots,t_n)\equiv1$, then we obtain what we call the {\it proper} non-Gaussian integral, which we denote simply by 
\begin{equation}
I(P)=I_{\alpha\equiv 1}(P)=\oint e^{P(t_1,\ldots,t_n)}dt_1\ldots dt_n.
\end{equation}
Arbitrary integral (2) is expressed through $I(P)$ as
\begin{equation}
I_\alpha(P)=\frac1{(2\pi)^n}\oint I(P-\mu_1 t_1-\ldots-\mu_n t_n)\widehat\alpha(\mu_1,\ldots,\mu_n)d\mu_1\ldots d\mu_n,
\end{equation}
where $\widehat\alpha(\mu_1,\ldots,\mu_n)$ is the Fourier--Laplace transform of $\alpha(t_1,\ldots,t_n)$. Function
$I(P-\mu_1 t_1-\ldots-\mu_n t_n)$ for fixed $P$ was considered in [2].

{\it Example} 1.4. Let $\alpha(t_1,\ldots,t_n)=t_1^{u_1-1}\ldots t_n^{u_n-1}$, where $u_1,\ldots,u_n$ are complex numbers. Then integral (2)
 has been studied in [3] and called a {\it $GG$-function}. We denote it by
\begin{equation}
I_{t_1^{u_1-1}\ldots t_n^{u_n-1}}(P)=\GG(P;u_1,\ldots,u_n)=\GG(c_\omega;u_1,\ldots,u_n).
\end{equation}
Arbitrary integral (2) is expressed through the $GG$-function as
\begin{equation}
I_\alpha(P)=\oint\GG(P;u_1,\ldots,u_n)\widetilde\alpha(u_1,\ldots,u_n)du_1\ldots du_n,
\end{equation}
where $\widetilde\alpha(u_1,\ldots,u_n)$ is the inverse Mellin transform of $\alpha(t_1$, $\ldots$, $t_n)t_1$ $\ldots$ $t_n$.

{\it Example} 1.5. Let $P_1,\ldots,P_k$ be polynomials in $t_1,\ldots,t_n$, let $\lambda_1$, $\ldots$, $\lambda_k$ be additional variables, let $P(\lambda_1,\ldots,\lambda_k,t_1,\ldots,t_n)$ be defined by the Cayley trick (5), 
and let $\alpha(\lambda_1,\ldots,\lambda_k,t_1,\ldots,t_n)$ be a function. Then we have 
\begin{equation}
\begin{aligned}
I_\alpha(P)&=\oint e^{\lambda_1 P_1+\ldots+\lambda_k P_k}\alpha(\lambda_1,\ldots,\lambda_k,t_1,\ldots,t_n)d\lambda_1\ldots d\lambda_k dt_1\ldots dt_n\\
&=\oint\beta(P_1(t_1,\ldots,t_n),\ldots,P_k(t_1,\ldots,t_n),t_1,\ldots,t_n)dt_1\ldots dt_n,
\end{aligned}
\end{equation}
where $\beta(y_1,\ldots,y_k,t_1,\ldots,t_n)$ is the Fourier--Laplace transform of $\alpha(\lambda_1$, $\ldots$, $\lambda_k,t_1,\ldots,t_n)$ with respect to $\lambda_1$, $\ldots$, $\lambda_k$.

{\bf Definition.} We call integral (11) by a {\it \emph(general\emph)  hypergeometric function} of coefficients of the polynomials $P_1,\ldots,P_k$, and denote it by 
$I_\beta(P_1,\ldots,P_k)$.

\section{Examples of general  hypergeometric functions}

Recall that by a general hypergeometric function we mean integral (3) or (11). 

{\it Example} 2.1. Let $k=1$ and $\beta(y,t_1,\ldots,t_n)=e^y\gamma(t_1,\ldots,t_n)$, where $\gamma(t_1,\ldots,t_n)$ is any function, then we obtain non-Gaussian integral (2) (with $\gamma$ instead of $\alpha$). 

{\it Example} 2.2. (Cf. [7, 8])  Let $y_{01},\ldots,y_{0k}$ be real numbers, let 
\begin{equation}
\beta(y_1,\ldots,y_k,t_1,\ldots,t_n)=\theta(y_{01}-y_1)\ldots\theta(y_{0k}-y_k)\gamma(t_1,\ldots,t_n),
\end{equation}
where $\gamma(t_1,\ldots,t_n)$ is any function,
\begin{equation}
\theta(y)=0\text{ if } y<0\text{ and }\theta(y)=1\text{ if }y\ge 0,
\end{equation}
and let $P_1,\ldots,P_k$ be polynomials with real coefficients. 
Then we obtain that the integral of $\gamma(t_1,\ldots,t_n)$ over the semi-algebraic domain 
\begin{equation}
P_1(t_1,\ldots,t_n)\le y_{01},\ \ \ldots,\ \ P_k(t_1,\ldots,t_n)\le y_{0k}
\end{equation}  
in the $n$-dimensional real space is a  hypergeometric function of coefficients of the polynomials $P_1,\ldots,P_k$.

{\it Example \emph{2.3:} the generalized Sturmfels theorem.} (Cf. [13, 12]) Let $k=1$ and
\begin{equation}
\beta(y,t_1,\ldots,t_n)=-\frac1{2\pi i}\log(y-y_0)\frac{\partial\gamma}{\partial t_1}(t_1,\ldots,t_n), 
\end{equation}
where $\gamma(t_1,\ldots,t_n)$ is any function.
Integrating by parts, we obtain that the integral
\begin{equation}
\begin{aligned}
-\frac1{2\pi i}&\oint\log(P(t_1,\ldots,t_n)-y_0)\frac{\partial\gamma}{\partial t_1}(t_1,\ldots,t_n)dt_1\ldots dt_n\\
=&\frac1{2\pi i}\oint\frac{\partial P/\partial t_1}{P(t_1,\ldots,t_n)-y_0}\gamma(t_1,\ldots,t_n)dt_1\ldots dt_n\\
=&\oint_{P(t_1,\ldots,t_n)=y_0}\gamma(t_1,\ldots,t_n)dt_2\ldots dt_n
\end{aligned}
\end{equation}
is a  hypergeometric function of coefficients of the polynomial $P$. In particular, if $n=1$, then we obtain the following theorem.
\medskip
 
{\bf Theorem 2.1.} {\it For a root $x$ of a polynomial equation $P(t)=y_0$ and for any function $\gamma(t)$, the quantity $\gamma(x)$ is a \emph(multi-valued\emph) hypergeometric function of coefficients of the polynomial $P$.}
\medskip

{\it Example} 2.4. (Cf. [13, 12]) Let $k=1$ and
\begin{equation}
\begin{aligned}
&\beta(y,t_1,\ldots,t_n)=\delta(y-y_0)\gamma(t_1,\ldots,t_n)\text{ or }\\
&\beta(y,t_1,\ldots,t_n)=\frac1{2\pi i(y-y_0)}\gamma(t_1,\ldots,t_n), 
\end{aligned}
\end{equation}
where $\gamma(t_1,\ldots,t_n)$ is any function. Then we obtain that the Gelfand--Leray integral [6] 
\begin{equation}
\oint_{P(t_1,\ldots,t_n)=y_0}\gamma(t_1,\ldots,t_n)dt_1\ldots dt_n/dP
\end{equation}
is a  hypergeometric function of coefficients of the polynomial $P$. 
Another proof of this fact follows from the generalized Sturmfels theorem: it suffices to differentiate equality (16) with respect to $y_0$ and replace $\partial\gamma/\partial t_1$ with $\gamma$. 

In particular, if $n=1$, then we obtain the following theorem.
\medskip

{\bf Theorem 2.2.} {\it For a root $x$ of a polynomial equation $P(t)=y_0$ and for any function $\gamma(t)$, the quantity $\gamma(x)/P'(x)$ is a \emph(multi-valued\emph)
hypergeometric function of coefficients of the polynomial $P$.}
\medskip

{\it Example} 2.5. (Cf. [5]) Let 
\begin{equation}
\beta(y_1,\ldots,y_k,t_1,\ldots,t_n)=y_1^{v_1}\ldots y_k^{v_k}t_1^{u_1-1}\ldots t_n^{u_n-1}.
\end{equation}
Then we obtain that the {\it generalized Euler integral} 
\begin{equation}
\oint P_1(t_1,\ldots,t_n)^{v_1}\ldots P_k(t_1,\ldots,t_n)^{v_k}t_1^{u_1-1}\ldots t_n^{u_n-1}dt_1\ldots dt_n
\end{equation}
is a  hypergeometric function of coefficients of the polynomials $P_1,\ldots,P_k$.

{\it Example} 2.6. (Cf. [13, 12]) Let $l\le k$ and
\begin{equation}
\beta(y_1,\ldots,y_k,t_1,\ldots,t_n)=\delta(y_1)\ldots\delta(y_l)\gamma(y_{l+1},\ldots,y_k,t_1,\ldots,t_n), 
\end{equation}
where $\gamma(y_{l+1},\ldots,y_k,t_1,\ldots,t_n)$ is any function. Then we obtain that
the Gelfand--Leray integral [6]
\begin{equation}
\oint\limits_{P_1=\ldots=P_l=0}\gamma(P_{l+1},\ldots,P_k,t_1,\ldots,t_n)dt_1\ldots dt_n/dP_1\ldots dP_l
\end{equation} 
is a  hypergeometric function of coefficients of the polynomials $P_1,\ldots,P_k$. In particular, if $l=n$, then we obtain that 
for a solution $(x_1$, $\ldots$, $x_n)$ of the system of equations 
\begin{equation}
P_1(x_1,\ldots,x_n)=\ldots=P_n(x_1,\ldots,x_n)=0,
\end{equation}
the quantity
\begin{equation}
\gamma(P_{n+1}(x_1,\ldots,x_n),\ldots,P_k(x_1,\ldots,x_n),x_1,\dots,x_n)/J(x_1,\ldots,x_n)
\end{equation}
is a  hypergeometric function of coefficients of the polynomials $P_1,\ldots,P_k$, where
\begin{equation}
J(x_1,\ldots,x_n)=\det(\partial P_i/\partial t_j)_{1\le i,j\le n}(x_1,\ldots,x_n)
\end{equation}
is the Jacobian of the polynomials $P_1,\ldots,P_n$ at the point  $t_1=x_1$, $\ldots$, $t_n=x_n$.
\medskip

\section{Equations satisfied by non-Gaussian integrals and by general hypergeometric functions}

\subsection{Equations satisfied by the non-Gaussian integral} 

{\bf Proposition 3.1.} {\it Non-Gaussian integral \emph{(2)} satisfies the following system of equations\emph: }
\begin{equation}
\frac{\partial I_\alpha}{\partial c_\omega}(P)=I_{t^\omega\alpha}(P) 
\end{equation}
{\it for any $\omega$, provided that the integral is regular in $c_\omega$, i.~e. admits differentiation under the sign of integral\emph;} 
\begin{equation}
I_{\frac{\partial\alpha}{\partial t_j}}(P)=-I_{\alpha\frac{\partial P}{\partial t_j}}(P),\ \ j=1,\ldots,n.
\end{equation}

{\bf Corollary.} {\it Non-Gaussian integral \emph{(2)} satisfies the equations }
\begin{equation}
\frac{\partial I_\alpha}{\partial c_\omega}=\frac{\partial^{\omega^1}}{\partial c_1^{\omega^1}}\ldots\frac{\partial^{\omega^n}}{\partial c_n^{\omega^n}}I_\alpha
\end{equation}
{\it for any $\omega$, where $c_j$ is the coefficient before the linear monomial $t_j$ in $P$}, $j=1,\ldots,n$;
\begin{equation}
I_{t_j\frac{\partial\alpha}{\partial t_j}}=-I_\alpha-\sum\limits_\omega\omega^j c_\omega \frac{\partial  I_\alpha}{\partial c_\omega}, \ \ j=1,\ldots,n;
\end{equation}
\begin{equation}
\frac{\partial}{\partial c_{\omega_1}}\ldots\frac{\partial}{\partial c_{\omega_N}}I_\alpha=\frac{\partial}{\partial c_{\omega'_1}}\ldots\frac{\partial}{\partial c_{\omega'_{N'}}}I_\alpha
\end{equation}
{\it for any $N, N'$ and any $\omega_1,\ldots,\omega_N,\omega'_1,\ldots,\omega'_{N'}$ such that }
\begin{equation}
\omega_1+\ldots+\omega_N=\omega'_1+\ldots+\omega'_{N'}.
\end{equation}
\medskip

System (26, 29) almost coincides with the $\GG$-system from [3]. 
System (29, 30) almost coincides with a corollary of the $\GG$-system called in [3] by the $A$-system with 
variables $c_\omega$, $\omega\in A$, where $A$ is a finite set of exponents $\omega=(\omega^1,\ldots,\omega^n)\in\Z^n$
of monomials $t^\omega=t_1^{\omega^1}\ldots t_n^{\omega^n}$. The $A$-system consists of equations (30) and the equations
\begin{equation}
\sum\limits_\omega\omega^j c_\omega\frac{\partial I_\alpha}{\partial c_\omega}=-u_j I_\alpha, \ \ j=1,\ldots,n,
\end{equation}
where $u_j$ are complex numbers (parameters). The $\GG$-system and the $A$-system are satisfied by the $\GG$-function (9).

\subsection{Equations satisfied by the general  hypergeometric function}

{\bf Proposition 3.2.} {\it General  hypergeometric integral \emph{(3, 11)}, considered as a function of coefficients $c_\omega^{(i)}$ of polynomials $P_i$, $i=1,\ldots,k$, 
satisfies the following system of equations\emph: }
\begin{equation}
\frac{\partial I_\beta}{\partial c_\omega^{(i)}}(P_1,\ldots,P_k)=I_{t^\omega\frac{\partial\beta}{\partial y_i}}(P_1,\ldots,P_k)
\end{equation}
{\it for any $\omega$ and $i$, provided that the integral is regular in $c_\omega^{(i)}$, i.~e. admits differentiation under the sign of integral\emph;} 
\begin{equation}
I_{y_i\beta}(P_1,\ldots,P_k)=I_{P_i\beta}(P_1,\ldots,P_k),\ \ i=1,\ldots,k;
\end{equation}
\begin{equation}
I_{\frac{\partial\beta}{\partial t_j}}(P_1,\ldots,P_k)=-\sum\limits_{i=1}^k I_{\frac{\partial\beta}{\partial y_i}\frac{\partial P_i}{\partial t_j}}(P_1,\ldots,P_k),\ \ j=1,\ldots,n.
\end{equation}

{\bf Corollary.} {\it General  hypergeometric integral satisfies the equations }
\begin{equation}
\left(\frac{\partial}{\partial c_0^{(i)}}\right)^{\omega^1+\ldots+\omega^n-1}\frac{\partial}{\partial c_\omega^{(i)}}I_\beta
=\left(\frac{\partial}{\partial c_1^{(i)}}\right)^{\omega^1}\ldots\left(\frac{\partial}{\partial c_n^{(i)}}\right)^{\omega^n}I_\beta
\end{equation}
{\it for any $\omega$ and $i$, where $c_0^{(i)}$ is the constant term of $P_i$ and $c_j^{(i)}$ is the coefficient before the linear monomial $t_j$ in $P_i$}, $i=1,\ldots,k$, $j=1,\ldots,n$;
\begin{equation}
I_{y_i\frac{\partial\beta}{\partial y_i}}=\sum\limits_\omega c_\omega^{(i)}\frac{\partial I_\beta}{\partial c_\omega^{(i)}},\ \ i=1,\ldots,k;
\end{equation}
\begin{equation}
I_{t_j\frac{\partial\beta}{\partial t_j}}=-I_\beta-\sum\limits_{\omega,i}\omega^j c_\omega^{(i)} \frac{\partial  I_\beta}{\partial c_\omega^{(i)}}, \ \ j=1,\ldots,n;
\end{equation}
\begin{equation}
\frac{\partial}{\partial c_{\omega_1}^{(i_1)}}\ldots\frac{\partial}{\partial c_{\omega_N}^{(i_N)}}I_\beta=\frac{\partial}{\partial c_{\omega'_1}^{(i_1)}}\ldots\frac{\partial}{\partial c_{\omega'_N}^{(i_N)}}I_\beta
\end{equation}
{\it for any $N$, any $i_1,\ldots,i_N$ and any $\omega_1,\ldots,\omega_N,\omega'_1,\ldots,\omega'_N$ such that }
\begin{equation}
\omega_1+\ldots+\omega_N=\omega'_1+\ldots+\omega'_N.
\end{equation}
\medskip

System (37--39) almost coincides with the $\widetilde A$-hypergeometric system ($\widetilde A$-system) from [3--5] with 
variables $c_\omega^{(i)}$, $\omega\in A_i$, $i=1,\ldots,k$, where $A_i$ is a finite set of exponents $\omega\in\Z^n$, and 
the set $\widetilde A\subset\Z^{n+k}=\Z^n\times\Z^k$ is defined as 
\begin{equation}
\widetilde A=A_1\times\{e_1\}\cup\ldots\cup A_k\times\{e_k\},
\end{equation}
where $e_1,\ldots,e_k$ is the standard basis in $\Z^k$ (the Cayley trick [5]). The $\widetilde A$-system consists of equations (39) and the equations
\begin{equation}
\sum\limits_\omega c_\omega^{(i)}\frac{\partial I_\beta}{\partial c_\omega^{(i)}}=v_i I_\beta,\ \ i=1,\ldots,k;
\end{equation}
\begin{equation}
\sum\limits_{\omega,i}\omega^j c_\omega^{(i)} \frac{\partial  I_\beta}{\partial c_\omega^{(i)}}=-u_jI_\beta, \ \ j=1,\ldots,n,
\end{equation}
where $v_i, u_j$ are complex numbers (parameters). The $\widetilde A$-system is satisfied by the generalized Euler integral (20).

\section{Power series expansions}

\subsection{Power series expansions of non-Gaussian integrals}

Let
\begin{equation}
P_0(t_1,\ldots,t_n)=\sum\limits_{\omega\in A} c_\omega^0 t^\omega,
\end{equation}
\begin{equation}
P(t_1,\ldots,t_n)=\sum\limits_{\omega\in A}c_\omega t^\omega=P_0(t_1,\ldots,t_n)+\sum\limits_{\omega\in A}a_\omega t^\omega,
\end{equation}
where $A\subset\Z^n$ is a finite set. Assume that non-Gaussian integral $I_\alpha(P)$ (2) is regular in a neighborhood of $P_0$. 

Following [3], let us call a set of exponents $B=\{\omega_1,\ldots,\omega_n\}\subset  A$ a {\it base} if they are linearly independent, i.~e. if they form a basis in $\C^n$.

We shall give the expansion of $I_\alpha(P)=I_\alpha(c_\omega)_{\omega\in A}$ into a power series in the variables
\begin{equation}
a_\omega=c_\omega-c_\omega^0,\ \ \omega\in A\setminus B,
\end{equation} 
with coefficients being functions of  $a_j=a_{\omega_j}=c_{\omega_j}-c_{\omega_j}^0$, $j=1,\ldots,n$. 

To this end, let us make the change of variables
\begin{equation}
T_j=t^{\omega_j},\ \ j=1,\ldots,n,
\end{equation}
in integral (2). We obtain 
\begin{equation}
\begin{aligned}
I_\alpha(P)&=\oint e^{\sum\limits_{j=1}^n a_j t^{\omega_j}+\sum\limits_{\omega\in A\setminus B}a_\omega t^\omega}\alpha_1(t_1,\ldots,t_n)dt_1\ldots dt_n\\
&=\oint e^{\sum\limits_{j=1}^n a_j T_j+\sum\limits_{\omega\in A\setminus B}a_\omega T^{l_\omega}}\alpha_2(T_1,\ldots,T_n)dT_1\ldots dT_n,
\end{aligned}
\end{equation}
where 
\begin{equation}
\alpha_1(t_1,\ldots,t_n)=e^{P_0(t_1,\ldots,t_n)}\alpha(t_1,\ldots,t_n),
\end{equation}
\begin{equation}
\alpha_2(T_1,\ldots,T_n)dT_1\ldots dT_n=\alpha_1(t_1,\ldots,t_n)dt_1\ldots dt_n, 
\end{equation}
and $l_\omega=(l_\omega^1,\ldots,l_\omega^n)$ is the vector of coordinates of $\omega$ with respect to the basis $\omega_1,\ldots,\omega_n$, 
\begin{equation}
\sum\limits_{j=1}^n l_\omega^j \omega_j=\omega.
\end{equation}
The numbers $l_\omega^j$ are, in general, rational numbers.

Expanding (48) into a power series in $a_\omega$, we obtain
\begin{equation}
I_\alpha(P)=\sum\limits_{\begin{subarray}{c}m_\omega\ge0\\ \omega\in A\setminus B\end{subarray}}C_m(a_1,\ldots,a_n)\prod\limits_{\omega}\frac{a_\omega^{m_\omega}}{m_\omega!},
\end{equation}
where
\begin{equation}
C_m(a_1,\ldots,a_n)=\oint e^{\sum\limits_{j=1}^n a_jT_j}T^{\sum\limits_{\omega\in A\setminus B}m_\omega l_\omega}\alpha_2(T_1,\ldots,T_n)dT_1\ldots dT_n
\end{equation}
is the Fourier--Laplace transform of $T^{\sum\limits_{\omega\in A\setminus B}m_\omega l_\omega}\alpha_2(T_1,\ldots,T_n)$.

In particular, if $\alpha(t_1,\ldots,t_n)=t_1^{u_1-1}\ldots t_n^{u_n-1}$ and $P_0=0$, then we obtain the expansion of $\GG$-function (9) into a power series of hypergeometric type [3].

If $B$ is the standard basis in $\Z^n$, then we obtain the obvious expansion of $I_\alpha(P)$ into a power series in $a_\omega$, $\omega\in A$,
\begin{equation}
I_\alpha(P)=\sum\limits_{\begin{subarray}{c}m_\omega\ge 0\\ \omega\in A\end{subarray}}\prod\limits_{\omega}\frac{a_\omega^{m_\omega}}{m_\omega!}I_{t^{\sum m_\omega\omega}\alpha}(P_0).
\end{equation}

\subsection{Power series expansions of general hypergeometric functions}  

According to formulas (4, 11), general hypergeometric integral $I_\beta(P_1,\ldots,P_k)$ (3) is a particular case of non-Gaussian integral $I_\alpha(P)$ (2,~5) corresponding to the set $\widetilde A\subset\Z^{n+k}$ defined by the 
Cayley trick (41). Hence, applying to this particular case the argument from \S4.1, we obtain, for any base $B\subset\widetilde A$, a power series expansion
of $I_\beta(P_1,\ldots,P_k)$. In particular, if $\beta(y_1,\ldots,y_k,t_1,\ldots,t_n)$ is given by (19), then this yields expansions of generalized Euler integral (20) into power series of hypergeometric type [4,~5].

\end{document}